\documentclass[12pt]{amsart}
\usepackage[cp866]{inputenc}
\usepackage{lcy}
\usepackage{amssymb}
\usepackage{verbatim}
\renewcommand{\leq}{\leqslant}
\renewcommand{\geq}{\geqslant}

\begin{document}

\title[Extending representations of braid groups]{Extending representations of braid groups to
the automorphism groups of free groups}

\author{Valerij G. Bardakov}

\address{Valerij
Bardakov, Sobolev Institute of Mathematics,
Novosibirsk,
630090, Russia}

\email{bardakov@math.nsc.ru}
\maketitle

\begin{abstract}
We construct a linear representation of the group  ${\rm IA}(F_n)$
of IA-automorphisms of a free group $F_n$, an extension of the
Gassner representation of the pure braid group $P_n$. Although
the problem of faithfulness of the Gassner representation is
still open for $n > 3$, we prove that the restriction of our
representation to the group of basis conjugating automorphisms
$Cb_n$ contains a non-trivial kernel even if $n = 2$.  We
construct also an extensions of the Burau representation to the
group of conjugating automorphisms $C_n$.  This representation is
not faithful for $n \geq 2$.
\end{abstract}


One of the generalizations
of the  braid group $B_n$ on the $n$ strings is the
group of conjugating automorphisms $C_n$.
The pure braid group $P_n$ is a normal subgroup of
the group $B_n.$
Similarly,  the group of basis
conjugating automorphisms $Cb_n$ is normal in
the group $C_n.$
In both cases, the quotient groups $B_n/P_n$
and $C_n/Cb_n$ are isomorphic to the
group $S_n,$ the symmetric group of degree $n.$

A.~G.~Savuschkina [1] proved that $C_n$ is a
semidirect product: $C_n = Cb_n \leftthreetimes S_n$.
For the group $B_n$ the similar statement is not
true, since $B_n$ is torsion-free.  The group of basis
conjugating automorphisms $Cb_n$ is a subgroup of the
group  $\operatorname{IA}(F_n)$ of the IA-automorphisms of a free
group $F_n$.

Naturally, the solution of the problem of the
linearity of the braid groups $B_n$ for all $n \geq 2$
[2, 3] initiated the study of the problem of the
linearity of $C_n$ (as well as an equivalent
problem of linearity of $Cb_n$) [4, Problem 15.9] and
of a more general problem of the linearity of the
group $\operatorname{IA}(F_n)$ [4,
Problem 15.10].

One of the possible approaches
to solution of these problems is to try
to extend
known linear representations of $P_n$ on $Cb_n$
(on $\operatorname{IA}(F_n)$) and to try to extend known
representations of $B_n$ on $C_n$. The most famous linear
representation of $B_n$ is the Burau representation [5] which is
faithful for $n=3$ and has a non-trivial kernel for all $n > 4$
[6, 7, 8]. In the case when $n=4,$ it is not known whether the Burau
representation of $B_n$ is faithful. There is a linear
representation of $P_n$ which is called the {Gassner}
representation. The problem of the faithfulness
of this representation for $n > 3$ is still open. Construction
of both these representations
stems from a general construction
of the so-called Magnus representation [9, Ch. 3].

In this article  using the Magnus representation we
will construct a linear representation
$\operatorname{IA}(F_n) \longrightarrow {\rm
GL}_n(\mathbb{Z}[t_1^{\pm 1}, t_2^{\pm 1},$ $\ldots, $
$t_n^{\pm 1}])$, which is an extension  of the Gassner
representation of $P_n$.  Then we can easily obtain
an extension of the Burau representation on $C_n$.
Unfortunately, the latter representation is not faithful.
Moreover, we will show that its restriction to the
group of basis-conjugating automorphisms $Cb_n$
has a non-trivial kernel even if $n = 2$, and so
the extension of the Burau representation to $C_n$
is not faithful for $n \geq 2$.

As we noted above the Burau representation is not
faithful for $n > 4$. However, we apply it to
calculate the Alexander polynomials of knots which are
closures of the corresponding braids.
The Alexander polynomials are one of the most
important invariants of the knots.
Similarly, $C_n$ is closely related to the  so-called
welded knots and links, and the linear
representations of $Cb_n$ and $C_n$ we have
constructed might be used to determine invariants of
the welded links.

{\bf {\it Acknowledgements.}} I would like to thank Vladimir Tolstykh for his remarks
on the first draft of this paper. Special thanks goes to the participants of the seminar
``Evariste Galois'' at Novosibirsk State University for their
kind attention to my work.

\section{Preliminaries}

The braid group $B_n$, $n\geq 2$, with $n$ strings can be
defined as the group generated by $n-1$ elements
$\sigma_1,\sigma_2,...,\sigma_{n-1}$ with the
following defining relations
\begin{align}
& \sigma_i\sigma_{i+1}\sigma_i=\sigma_{i+1}\sigma_i\sigma_{i+1},~~~ i=1,2,...,n-2, \\
& \sigma_i \sigma_j = \sigma_j \sigma_i,~~~ |i-j|\geq 2.
\end{align}

There is a homomorphism from the group $B_n$ to the
symmetric group $S_n$ of degree $n$ defined via
$$
\nu(\sigma_i)=(i,i+1),~~~i=1,2,\ldots,n-1.
$$
The kernel of the homomorphism $\nu $ is the pure braid group $P_n$.
The group $P_n$ admits a presentation with generators
\begin{align*}
& a_{i,i+1}=\sigma_i^2, \\
& a_{ij}=\sigma_{j-1}\sigma_{j-2}\ldots\sigma_{i+1}\sigma_i^2\sigma_{i+1}^{-1}\ldots
\sigma_{j-2}^{-1}\sigma_{j-1}^{-1},~~~i+1< j \leq n,
\end{align*}
and the following defining relations:
$$
\begin{array}{ll}
a_{ik}^{-\varepsilon }a_{kj}a_{ik}^{\varepsilon }=\left( a_{ij}a_{kj}\right) ^{\varepsilon}a_{kj}
\left( a_{ij}a_{kj}\right) ^{-\varepsilon }, & \\
& \\
a_{km}^{-\varepsilon }a_{kj}a_{km}^{\varepsilon }=\left( a_{kj}a_{mj}\right) ^{\varepsilon }a_{kj}
\left( a_{kj}a_{mj}\right) ^{-\varepsilon }, &~~~m < j,\\
& \\
a_{im}^{-\varepsilon }a_{kj}a_{im}^{\varepsilon }=\left[ a_{ij}^{-\varepsilon }, a_{mj}^{-\varepsilon }\right] ^{\varepsilon }a_{kj}
\left[ a_{ij}^{-\varepsilon }, a_{mj}^{-\varepsilon}\right] ^{-\varepsilon}, & ~~~i < k <
m,\\
& \\
a_{im}^{-\varepsilon }a_{kj}a_{im}^{\varepsilon }=a_{kj}, & ~~~k < i;~~~m < j~~
\mbox{or}~~ m < k,
\end{array}
$$
where $\varepsilon = \pm 1$.

The subgroup $U_n$ generated by $a_{1,n},
a_{2,n},\ldots,a_{n-1,n}$ is a free normal subgroup
of $P_n$. The group $P_n$ is a
semidirect product of $U_n$ and $P_{n-1}$. Hence the group $P_n$
is a semi direct product
$$
P_n=U_n\leftthreetimes (U_{n-1}\leftthreetimes (\ldots \leftthreetimes
(U_3\leftthreetimes U_2))\ldots),
$$
where $U_i=\langle a_{1,i},
a_{2,i},\ldots,a_{i-1,i} \rangle, i=2,3,\ldots,n,$ is a free group of
rank $i-1$.

The braid group $B_n$ can be embedded into
the automorphism group of a free group $F_n$
with a free basis $\{x_1, x_2,\ldots, x_n \}.$
The said embedding is induced by a map from $B_n$
to ${\rm Aut}(F_n)$ defined by
$$
\sigma_{i} : \left\{
\begin{array}{ll}
x_{i} \longmapsto x_{i}x_{i+1}x_i^{-1}, &  \\ x_{i+1} \longmapsto
x_{i}, & \\ x_{l} \longmapsto x_{l} & ~~ l\neq i,i+1.
\end{array} \right.
$$
The generator $a_{rs}$ of $P_n$ determines the following automorphism
of $F_n$
$$
a_{rs} : \left\{
\begin{array}{ll}
x_{i} \longmapsto x_{i} & \mbox{if }~~ s < i~~ \mbox{or}~~ i < r, \\
& \\
x_{r} \longmapsto x_{r}x_{s}x_rx_s^{-1}x_r^{-1}, &  \\
& \\
x_{i} \longmapsto [x_{r}^{-1}, x_s^{-1}]x_i[x_{r}^{-1}, x_s^{-1}]^{-1} & \mbox{if }~~ r < i <s,\\
& \\
x_{s} \longmapsto x_{r}x_sx_r^{-1}. &
\end{array} \right.
$$

By a theorem of Artin [9, Theorem 1.9] automorphism $\beta$
in ${\rm Aut}(F_n)$ belongs to (the image of) $B_n$
if and only if $\beta$ satisfies the following
two conditions
\begin{alignat*}3
&1)\quad &&\beta(x_i)=f_i^{-1}x_{\pi(i)}f_i,            && 1\leq i\leq n,\\
&2) &&\beta(x_1x_2 \ldots x_n)=x_1x_2 \ldots x_n,  &&
\end{alignat*}
where $\pi $ is a permutation of  $\{1, 2, \ldots, n \}$ and $f_i=f_i(x_1,x_2,\ldots,x_n)$
is a word in the generators of $F_n$.

An automorphism of $F_n$ is called a {\it conjugating automorphism} if it
satisfies condition 1). Let $C_n$ be the group of conjugating
automorphisms. Evidently,
$Cb_n$ is normal in $C_n$ and  the quotient group
$C_n/Cb_n$ is isomorphic to the symmetric group $S_n$.
The elements of the group $Cb_n$ are called  {\it
basis-conjugating automorphisms}. J.~McCool [11]
proved that the group $Cb_n$ is generated by
automorphisms
$$
\varepsilon_{ij} : \left\{
\begin{array}{ll}
x_{i} \longmapsto x_{j}^{-1}x_ix_j, &  i\neq j, \\
x_{l} \longmapsto x_{l} &  l\neq i,
\end{array} \right.
$$
$i\leq i\neq j \leq n$.

Recall (see [12, chapter 1, \S~4]) that
$\operatorname{IA}(F_n)$ is generated by the
automorphisms $\varepsilon_{ij}$, $1\leq i\neq j \leq
n$ and the  automorphisms
$$
\varepsilon_{ijk} :
\begin{cases}
x_{i} \longmapsto x_{i}[x_j, x_k],  & \text{ if } k\neq i, j, \\
x_{l} \longmapsto x_{l},            & \text{ if } l\neq i,
\end{cases}
$$
where $[a, b]=a^{-1}b^{-1}ab.$

\vskip 20pt

\section{Fox's derivatives and Magnus representation}

Recall the  definitions and main properties of Fox's derivatives [9,
Chapter 3; 13, Chapter 7].

Let $F_n$ be a free group of rank $n$ with free generators $x_1, x_2, \ldots, x_n$.
If $\varphi$ is any homomorphism defined
on $F_n,$ then we use the symbol
$F_n^{\varphi}$ to denote the image of $F_n$ under $\varphi$.
Consider also the group ring $\mathbb{Z}F_n$ of
the group $F_n$ over
the ring $\mathbb{Z}$ of integers.

For every $j=1, 2, \ldots, n$ define the mapping
$$
\frac{\partial }{\partial x_j} : \mathbb{Z}F_n\longrightarrow \mathbb{Z}F_n
$$
using the following conditions
\begin{alignat*}2
& 1)\quad &&\frac{\partial x_i}{\partial x_j}=
\begin{cases}
1, & \text{ if } i = j , \\
0, & \text{ if } i\not= j,
\end{cases}
\\
& 2) && \frac{\partial x_i^{-1}}{\partial x_j}=
\begin{cases}
-x_i^{-1}, & \text{ if } i = j , \\
0,         & \text{ if } i\not= j,
\end{cases},\\
& 3) && \frac{\partial (w v)}{\partial x_j} = \frac{\partial w}{\partial x_j} (v)^{\tau } +
w \frac{\partial v}{\partial x_j},~~~w, v \in \mathbb{Z}F_n,
\end{alignat*}
where $\tau : \mathbb{Z}F_n \longrightarrow
\mathbb{Z}$ is the operation of trivialiazation which
sends all elements of $F_n$ to the identity,
$$
4)\quad \frac{\partial }{\partial x_j} \left( \sum a_g g\right) =
\sum a_g \frac{\partial g}{\partial x_j},~~~g\in F_n,~~~a_g \in
\mathbb{Z}.
$$

If we denote the fundamental ideal of  ring
$\mathbb{Z}F_n$ (the kernel of the homomorphism
$\tau$) by $\Delta_n$,
then it is easy to see that for
every $v \in \mathbb{Z}F_n$ the element $v-v^{\tau}$
belongs to $\Delta_n$.  The following formula is true:
$$
v-v^{\tau } = \sum\limits_{j=1}^n \frac{\partial v}{\partial x_j}
(x_j - 1);
$$
this formula is called the {\it ``Fundamental
formula'' of free calculus.} In particular, we
have as a consequence that $\{x_1 - 1, x_2 - 1, \ldots, x_n - 1\}$
is a basis of the fundamental ideal $\Delta_n$.

We shall need Blanchfield's theorem [9, Theorem
3.5] which says
that an element $v \in F_n$ lies in the commutator subgroup
$[\mbox{ker} \varphi, \mbox{ker} \varphi]$ if and only if
$\displaystyle  \left( \frac{\partial v}{\partial x_j}\right)^{\varphi} = 0$
for all  $j = 1, 2, \ldots, n$.

Let $A_{\varphi}$ be any subgroup of ${\rm Aut} F_n$ which satisfy the condition
$$
x^{\varphi} = (x^{\alpha})^{\varphi}
$$
for every $x\in F_n$ and for every $\alpha \in A_{\varphi}$. If $\alpha \in
A_{\varphi},$ we define $|| \alpha ||$ to be the $n \times m$
matrix
$$
|| \alpha || = \left[ \left( \frac{\partial (x_i^{\alpha })}{\partial
x_j} \right)^{\varphi} \right]_{i,j=1}^n
$$
with entries in
$\mathbb{Z}F_n^{\varphi}$. This mapping defines the Magnus representation
$$
\rho : A_{\varphi} \longrightarrow {\rm GL}_n
(\mathbb{Z}F_n^{\varphi}).
$$
Taking as $\varphi$ any homomorphism from $F_n$ onto
an infinite cyclic group $\langle t\rangle,$ that is,
assuming that
$$
x_i^{\varphi} = t,\quad i=1, 2, \ldots,n,
$$
we obtain the Burau representation
$$
\rho_B : B_{n} \longrightarrow {\rm GL}_n
(\mathbb{Z}[t^{\pm 1}])
$$
of the braid group $B_n$.

Similarly, if $\varphi$ is a homomorphism
from $F_n$ onto a free abelian group $A_n$
with free generators $t_1, t_2, \ldots, t_n$,
$$
x_i^{\varphi} = t_i, \quad i=1, 2, \ldots, n
$$
then the resulting representation is the Gassner representation
$$
\rho_G : P_{n} \longrightarrow {\rm GL}_n
(\mathbb{Z}[t_1^{\pm 1}, t_2^{\pm 1}, \ldots, t_n^{\pm 1}])
$$
of the group $P_n$.

\vskip 20pt

\section{The construction of $\widehat{\rho}_G$}

It is easy to check  that for every  $x \in F_n$ and every automorphism
$\alpha \in \operatorname{IA}(F_n)$  the following equality  is true:
$$
x^{\varphi} = (x^{\alpha})^{\varphi},
$$
where the homomorphism $\varphi : F_n \longrightarrow
A_n$ is defined in the end of the previous section.
Consequently, we can construct the Magnus
representation
$$
\rho : \operatorname{IA}(F_n) \longrightarrow {\rm GL}_n(R),
$$
where $R = \mathbb{Z}[t_1^{\pm 1}, t_2^{\pm 1},
\ldots, t_n^{\pm 1}]$.  In order to define the action
of $\rho $  on generators of $\operatorname{IA}(F_n)$
we should calculate the Fox's derivations.

{\bf Lemma 1.} {\it The following formulas are true in the group ring $\mathbb{Z}F_n$}
({\it where $i, j, k$ are pairwise distinct}):
$$
\begin{array}{ll}
\displaystyle \frac{\partial x_l^{\displaystyle \varepsilon_{ijk}}}{\partial x_m} = 0, & m \not\in \{l, i, j, k\},\\
& \\
\displaystyle \frac{\partial x_l^{\displaystyle \varepsilon_{ijk}}}{\partial x_l} = 1, &  l \not\in \{j, k\},\\
& \\
\displaystyle \frac{\partial x_i^{\displaystyle \varepsilon_{ijk}}}{\partial x_j} = -x_i x_j^{-1} + x_i x_j^{-1} x_k^{-1}, & \\
& \\
\displaystyle \frac{\partial x_i^{\displaystyle \varepsilon_{ijk}}}{\partial x_k} = -x_i x_j^{-1} x_k^{-1} + x_i x_j^{-1}
x_k^{-1} x_j. & \\
\end{array}
$$

We will  consider the matrix $\rho(\varepsilon_{ijk})$
as an automorphism of a free left $R$--module $W_n$
with base $e_1, e_2, \ldots, e_n$.

Then the action of this automorphism (from the right) on the
base vectors is as follows:
$$
\left\{
\begin{array}{ll}
e_{i} \rho(\varepsilon_{ijk}) = e_{i} + t_i t_j^{-1} (t_k^{-1} - 1)e_j +
 t_i t_k^{-1} (1 - t_j^{-1})e_k, &  \\
 & \\
e_{l} \rho(\varepsilon_{ijk}) = e_{l} & \mbox{if } l\neq i.
\end{array} \right.
$$

In order to construct the matrices
$\rho(\varepsilon_{ij})$, we will use the following
simple lemma.

{\bf Lemma 2.} {\it The following formulas are true in the group ring $\mathbb{Z}F_n$
}
({\it where $i, j, l, m$ are pairwise distinct}):
$$
~~~\displaystyle \frac{\partial x_l^{\displaystyle \varepsilon_{ij}}}{\partial x_m} =
0,
~~~\displaystyle \frac{\partial x_l^{\displaystyle \varepsilon_{ij}}}{\partial x_l} = 1,
~~~\displaystyle \frac{\partial x_i^{\displaystyle \varepsilon_{ij}}}{\partial x_i} = x_j^{-1},
~~~\displaystyle \frac{\partial x_i^{\displaystyle \varepsilon_{ij}}}{\partial x_j} = x_j^{-1}( x_i - 1).
$$
$$
\displaystyle \frac{\partial x_l^{\displaystyle \varepsilon_{ij}^{-1}}}{\partial x_m} = 0,
~~~\displaystyle \frac{\partial x_l^{\displaystyle \varepsilon_{ij}^{-1}}}{\partial x_l} = 1,
~~~\displaystyle \frac{\partial x_i^{\displaystyle \varepsilon_{ij}^{-1}}}{\partial x_i} = x_j,
~~~\displaystyle \frac{\partial x_i^{ \displaystyle \varepsilon_{ij}^{-1}}}{\partial x_j} = 1 - x_j x_i
x_j^{-1}.
$$

Then the matrix  $\rho(\varepsilon_{ij})$  defines the
action on the base of $W_n$ by the next formulas:
$$
\left\{
\begin{array}{ll}
e_{i} \rho(\varepsilon_{ij}) = t_j^{-1} (t_i - 1) e_{j} + t_j^{-1} e_i, &  \\
& \\
e_{l} \rho(\varepsilon_{ij}) = e_{l} & \mbox{if } l\neq i.
\end{array} \right.
$$

For the automorphism
$$
\varepsilon_{ij}^{-1} : \left\{
\begin{array}{ll}
x_{i} \longmapsto x_{j} x_i x_j^{-1} & \mbox{if } i\neq j, \\
 & \\
x_{l} \longmapsto x_{l} & \mbox{if } l\neq i,
\end{array} \right.
$$
the matrix $\rho(\varepsilon_{ij}^{-1})$  is defined by the action on the base of $W_n$
as follows:
$$
\left\{
\begin{array}{ll}
e_{i} \rho(\varepsilon_{ij}^{-1}) = t_j e_{i} + (1 - t_i)e_j, &  \\
 & \\
e_{l} \rho(\varepsilon_{ij}^{-1}) = e_{l} & \mbox{if } l\neq i.
\end{array} \right.
$$
It is easy to check that being defined in this way the
mapping $\rho$ is a linear representation of ${\rm
IA}(F_n)$.

Restricting $\rho$ to $Cb_n,$ we obtain that
$$
\widehat{\rho}_G = \rho \vert _{Cb_n} : Cb_n \longrightarrow {\rm
GL}_n (R).
$$

{\bf Theorem 1.} {\it The linear representation
$$
\widehat{\rho}_G : Cb_n \longrightarrow {\rm
GL}_n (R)
$$
is an extensions of the Gassner representation of the pure braid group $P_n$.
}

{\bf Proof.}
The group $P_n$ is a subgroup of $Cb_n$ and its generators $a_{ij}$,
$1 \leq i < j \leq n$ can be expressed
via the (standard) generators of $Cb_n$ in the following way [14, Lemma 4]:
\begin{alignat*}3
& a_{i,i+1} &&=\varepsilon_{i,i+1}^{-1}\varepsilon_{i+1,i}^{-1}, && i=1,2,\ldots,n-1, \\
& a_{ij} &&=\varepsilon_{j-1,i}\varepsilon_{j-2,i}\ldots\varepsilon_{i+1,i}
(\varepsilon_{ij}^{-1}\varepsilon_{ji}^{-1})\varepsilon_{i+1,i}^{-1}\ldots
\varepsilon_{j-2,i}^{-1}\varepsilon_{j-1,i}^{-1}  && \\
& &&=\varepsilon_{j-1,j}^{-1}\varepsilon_{j-2,j}^{-1}\ldots\varepsilon_{i+1,j}^{-1}
(\varepsilon_{ij}^{-1}\varepsilon_{ji}^{-1})\varepsilon_{i+1,j}\ldots
\varepsilon_{j-2,j}\varepsilon_{j-1,j}, \quad && 2 \leq i+1 < j \leq n.
\end{alignat*}

Using these formulas
we find the matrix  $\rho(a_{ij})$, $1 \leq i < j \leq n$
and then compare it with the  Gassner matrix from [9, p.118].

Let us first find the matrix
\begin{align*}
\rho(a_{i,i+1}) = \rho(\varepsilon_{i,i+1}^{-1}\varepsilon_{i+1,i}^{-1})
                =  \rho(\varepsilon_{i,i+1}^{-1}) \rho(\varepsilon_{i+1,i}^{-1}).
\end{align*}
We have
$$
\left\{
\begin{array}{ll}
e_{i} \rho(\varepsilon_{i,i+1}^{-1}) \rho(\varepsilon_{i+1,i}^{-1}) =
(t_{i+1} e_i + (1-t_i) e_{i+1}) \rho(\varepsilon_{i+1,i}^{-1}) = & \\
& \\
\qquad \qquad \qquad = (1 - t_i + t_i t_{i+1}) e_i +
 t_i (1 - t_i) e_{i+1}, &  \\
 & \\
 e_{i+1} \rho(\varepsilon_{i,i+1}^{-1}) \rho(\varepsilon_{i+1,i}^{-1}) =
e_{i+1} \rho(\varepsilon_{i+1,i}^{-1}) = (1 - t_{i+1}) e_i +
 t_i e_{i+1}, &  \\
 & \\
e_{l} \rho(\varepsilon_{i,i+1}^{-1}) \rho(\varepsilon_{i+1,i}^{-1}) = e_{l}
& \mbox{if}~~ l\neq i,
i+1.
\end{array} \right.
$$

Next, we have to find  $\rho(a_{ij})$ when $j > i+1$.
To do this we have to calculate the matrix
$\rho(\varepsilon_{i,j}^{-1} \varepsilon_{j,i}^{-1})$:
$$
\left\{
\begin{array}{ll}
e_{i} \rho(\varepsilon_{i,j}^{-1}) \rho(\varepsilon_{j,i}^{-1}) =
(t_{j} e_i + (1-t_i) e_{j}) \rho(\varepsilon_{j,i}^{-1}) = (1 - t_i + t_i t_{j}) e_i +
 t_i (1 - t_i) e_{j}, &  \\
 & \\
 e_{j} \rho(\varepsilon_{i,j}^{-1}) \rho(\varepsilon_{j,i}^{-1}) =
e_{j} \rho(\varepsilon_{j,i}^{-1}) = (1 - t_{j}) e_i +
 t_i e_{j}, &  \\
 & \\
e_{l} \rho(\varepsilon_{i,j}^{-1}) \rho(\varepsilon_{j,i}^{-1}) = e_{l} & \mbox{if } l\neq i, j.
\end{array} \right.
$$
The following equalities may  be easily to check:
$$
\left\{
\begin{array}{ll}
e_{l} \rho(a_{i,j}) = e_{l} & \mbox{if } l < i ~~~\mbox{or}~~~ l > j,
\\
& \\
e_{i} \rho(a_{i,j}) = (1 - t_i + t_i t_{j}) e_i +
 t_i (1 - t_i) e_{j}, &  \\
 & \\
e_{k} \rho(a_{i,j}) = (t_k - 1)(t_j - 1) e_i + e_k
+ (t_k - 1) (1 - t_i) e_{j} & \mbox{if } i < k < j, \\
& \\
e_{j} \rho(a_{i,j}) = (1 - t_{j}) e_i +
 t_i e_{j}. &  \\
\end{array} \right.
$$
Сomparing this matrix with the Gassner matrix,
we get the desired conclusion.

Recall that the group of conjugating automorphisms
$C_n$ can be decomposed
as a semidirect product $C_n
= Cb_n\leftthreetimes S_n$.

As the matrix $\widehat{\rho}_G(\varepsilon_{ij})$
depends on $t_1, t_2,\ldots, t_n,$ we may set
$$
t_1 = t_2 = \cdots = t_n = t,
$$
thereby obtaining the matrix which
we shall denote by  $\widehat{\rho}_B(\varepsilon_{ij})$.
Further, for each automorphism from  $S_n \leq \mbox{Aut}(F_n)$
we assign the matrix of the corresponding
permutation of the elements of the
base $W_n$. We will then obtain the representation
$$
\widehat{\rho}_B : C_n \longrightarrow {\rm
GL}_n(\mathbb{Z}[t^{\pm}]).
$$
It is easy to check that $\widehat{\rho}$ is an
extension of  the Burau representation of the braid group
$B_n$.

As $Cb_n$ is a subgroup of finite index in  $C_n$, the representation
$\widehat{\rho}_B$ is faithful if its restriction to
$Cb_n$ is faithful.

It is well-known that the Burau representation
$\rho_B$  and the Gassner representations $\rho_G$ are
reducible; however, both these representations
determine some irreducible representations
of dimension $n-1$ [9, Lemma 3.11.1].

The following question naturally arises.

{\bf Question.} Is it true that the presentations
$\widehat{\rho}_B$ and $\widehat{\rho}_G$ are reducible?

\section{Kernel of the representation $\widehat{\rho}_G$}

In this section  we will show that the representation $\widehat{\rho}_G$ is not faithful.
We shall prove the following result.

{\bf Theorem 2.} {\it The representation
$\widehat{\rho}_G : Cb_n \longrightarrow {\rm GL}_n(R)$
has a non-trivial kernel for every} $n \geq 2$.

The theorem implies that the representation $\rho :
\operatorname{IA}(F_n) \longrightarrow {\rm GL}_n(R)$
has a non-trivial kernel for every $n \geq 2$.

Recall [14] that the group of basis--conjugating
automorphisms  $Cb_n$, $n\geq 2$ can be decomposed
as a semidirect product as follows:
$$
Cb_n=D_{n-1}\leftthreetimes (D_{n-2}\leftthreetimes (\ldots
\leftthreetimes (D_2\leftthreetimes D_1))\ldots ),
$$
where subgroup $D_i$ is generated by the elements
$$
\varepsilon_{i+1,1},\varepsilon_{i+1,2},
\ldots,\varepsilon_{i+1,i}, \varepsilon_{1,i+1},
\varepsilon_{2,i+1},\ldots,\varepsilon_{i,i+1}.
$$
Moreover, the elements $\varepsilon_{i+1,1},\varepsilon_{i+1,2},
\ldots,\varepsilon_{i+1,i}$ generate a free group of  rank $i$ which we will  denote
by
$L_i$ and the elements $\varepsilon_{1,i+1},\varepsilon_{2,i+1},\ldots,\varepsilon_{i,i+1}$
generate a free abelian group of  rank $i$ which we will denote by $A_i$.

Let us show that if $j \geq 3$ then the second commutator subgroup $L_j''$
is contained in the kernel of the representation $\widehat{\rho}_G$.
Following the construction of the matrix
$\widehat{\rho}_G(\varepsilon_{ij})$ =
$\rho(\varepsilon_{ij}),$ one sees that
it different from the identity matrix only in $i$th row.
Let $w =w(\varepsilon_{i1}, \varepsilon_{i2}, \ldots,
\varepsilon_{i,i-1})$ be a reduced word
over the elements $\varepsilon_{i1}^{\pm 1},
\varepsilon_{i2}^{\pm 1}, \ldots,
\varepsilon_{i,i-1}^{\pm 1},$ the free
generators of $L_{i-1}$.  The following lemma describes
how the automorphism $w$ acts on a generator
$x_i$ of the free group $F_n$.

{\bf Lemma 3.} {\it Let $w = w(\varepsilon_{i1}, \varepsilon_{i2}, \ldots, \varepsilon_{i,i-1})$
be a reduced word which represents an element of $L_{i-1}$.
Then
$$
x_i w(\varepsilon_{i1}, \varepsilon_{i2}, \ldots, \varepsilon_{i,i-1}) =
\displaystyle x_i^{w^*(x_1, x_2, \ldots, x_{i-1})},
$$
where the word $w^*(x_1, x_2, \ldots, x_{i-1})$ is the
{\em reverse} word of the word $w(x_1, x_2, \ldots,
x_{i-1})$, that is the syllables of $w^*$ are
the syllables of $w$ written in the reverse order.}

{\bf Proof.}  The statement is a consequence
of the following equations that can be
easily verified by induction on the
number of syllables of $w$ (the {\it syllabic
length} of $w$):
$$
x_i \varepsilon_{ik}^p = x_i^{x_k^p},~~~x_i \varepsilon_{ik}^p \varepsilon_{il}^q =
\displaystyle x_i^{x_l^q x_{k}^p},~~~1 \leq k \neq l \leq i-1,~~~p, q \in \mathbb{Z},
$$

Obviously, if $w = w(\varepsilon_{i1},
\varepsilon_{i2}, \ldots, \varepsilon_{i,i-1})$
represents an element of the second commutator subgroup
$L_{i-1}'',$ then the word $w^* = w^*(x_1, x_2, \ldots,
x_{i-1})$ represents an element of the second
commutator subgroup $F_n''$.

Assume that
$w = w(\varepsilon_{i1}, \varepsilon_{i2}, \ldots, \varepsilon_{i,i-1})$
represents an element of $L_{i-1}''$. In order to find the matrix $\widehat{\rho}_G(w)$
we have to find
the derivatives $\displaystyle{\frac{\partial (x_i w)}{\partial x_k}}$, $k=1, 2, \ldots, n$.
Let us check the following equality:
$$
\left( \frac{\partial (x_i w)}{\partial x_k}\right) ^{\varphi} =
\left\{
\begin{array}{ll}
1, & \mbox{if}~~~ i = k,\\
& \\
0, &  \mbox{otherwise.}
\end{array}
\right.
$$
The case when $k > i$ and the case when the word $x_i w$ does not contain $x_k$
are simple. Suppose that $k=i.$ We have
$$
\frac{\partial (x_i w)}{\partial x_i} = \frac{\partial ((w^*)^{-1} x_i w^*)}{\partial
x_i}.
$$
Using properties of the Fox's derivatives, we obtain that
\begin{align*}
\frac{\partial ((w^*)^{-1} x_i w^*)}{\partial x_i} &=
\frac{\partial ((w^*)^{-1})}{\partial x_i} + (w^*)^{-1} \frac{\partial (x_i w^*)}{\partial x_i}
   \\
&=-(w^*)^{-1} \frac{\partial w^*}{\partial x_i} +
(w^*)^{-1} \left( 1 + x_i \frac{\partial w^*}{\partial x_i}\right).
\end{align*}
As $w^*$ does not contain $x_i,$ the derivatives
$\displaystyle{\frac{\partial w^*}{\partial x_i}}$
become zero and $w^*$ lies in the commutator subgroup
$F_n'$ and $\varphi$ takes it to 1,
that is, $(w^*)^{\varphi} = 1,$ as required.

Let, finally, $k < i$. Then
$$
\frac{\partial ((w^*)^{-1} x_i w^*)}{\partial x_k} =
(w^*)^{-1} (x_i - 1) \frac{\partial w^*}{\partial x_k}.
$$
Since ${\rm ker} \varphi = F_n',$ then by
Blanchfield's theorem (taking into account that $w^* \in
F_n''$) we have the equality
$\displaystyle{\left(\frac{\partial w^*}{\partial
x_k}\right)^{\varphi}} = 0$.  The
statement is proven.

As a consequence we have the following
fact.

{\bf Lemma 4.} {\it
The representation $\widehat{\rho}_G$ of $Cb_n$, $n
\geq 3$, is not faithful, since its kernel contains
the subgroups $L_i''$, $i=2, 3, \ldots, n-1$.  In
particular, the above constructed representation
$\widehat{\rho}_B$ of $C_n,$ an extension of
the Burau representation, is not faithful  for all} $n
\geq 3$.

Let us demonstrate that actually the  representation
$\widehat{\rho}_G$ is not faithful even if $n=2$.
In order to prove that, we shall need

{\bf Lemma 5.} {\it The following formulas
$$
x_i w(\varepsilon_{12}, \varepsilon_{21}) = x_i^{w(x_2, x_1)},~~~i
= 1, 2.
$$
are true in $F_2 = \langle x_1, x_2 \rangle $. In particular, if
$w(\varepsilon_{12}, \varepsilon_{21})$ is in the
second commutator subgroup $D_1'',$ then  $w(x_2, x_1)$
is in the second commutator subgroup} $F_2''$.

{\bf Proof.}
Recall that $Cb_2 = D_1 = \langle \varepsilon_{21}, \varepsilon_{12}\rangle \simeq F_2$
and the automorphisms
$\varepsilon_{12}$ and $\varepsilon_{21}$ are defined as
follows:
$$
\varepsilon_{12} : \left\{
\begin{array}{ll}
x_{1} \longmapsto x_{2}^{-1} x_1 x_2, &  \\
x_{2} \longmapsto x_{2}, &
\end{array} \right. ~~~~~
\varepsilon_{21}^{-1} : \left\{
\begin{array}{ll}
x_{1} \longmapsto x_{1}, &  \\
x_{2} \longmapsto x_{1}^{-1} x_2 x_1,  &
\end{array} \right.
$$
i. e.  $\varepsilon_{12}$ induces conjugation by $x_2$ in $F_2$ and
$\varepsilon_{21}$ conjugation by $x_1$.
The action of $\varepsilon_{12}\varepsilon_{21}$
on the base $x_1,x_2$ is
$$
\varepsilon_{12}\varepsilon_{21} : \left\{
\begin{array}{ll}
x_{1} \longmapsto x_{1}^{-1} x_2^{-1} x_1 x_2 x_1, &  \\
x_{2} \longmapsto x_{1}^{-1} x_2 x_1,  &
\end{array} \right.
$$
which means that $x_i \varepsilon_{12} \varepsilon_{21} = x_i^{x_2 x_1}$, $i=1,2$.

If $w = w(\varepsilon_{12}, \varepsilon_{21})$ is a
reduced word over the alphabet $\{ \varepsilon_{12}^{\pm
1}, \varepsilon_{21}^{\pm 1} \},$ then using induction
on the length of $w,$ we get the following
statement.

{\bf Lemma 6.} {\it The kernel of the representation $\widehat{\rho}_G$ of
$Cb_2$
 coincides the second commutator subgroup} $D_1''$.

{\bf Proof.}
Let us show that if $w = w(\varepsilon_{12}, \varepsilon_{21})$ is in $D_1''$
then
$\widehat{\rho}_G(w) = 1$.
It follows from the definition of the Magnus representation
that
$$
\widehat{\rho}_G(w(\varepsilon_{12}, \varepsilon_{21})) =
\left(
\begin{array}{cc}
\displaystyle \frac{\partial (x_1 w(\varepsilon_{12}, \varepsilon_{21}))}{\partial x_1} &
\displaystyle \frac{\partial (x_1 w(\varepsilon_{12}, \varepsilon_{21}))}{\partial x_2}  \\
& \\
\displaystyle \frac{\partial (x_2 w(\varepsilon_{12}, \varepsilon_{21}))}{\partial x_1} &
\displaystyle \frac{\partial (x_2w(\varepsilon_{12}, \varepsilon_{21}))}{\partial x_2}   \\
\end{array}
\right) ^{\varphi}.
$$
Write $w_1$ for $w(x_2, x_1).$ We calculate
the Fox's derivatives:
$$
\frac{\partial (x_1 w)}{\partial x_1} = \frac{\partial (w_1^{-1} x_1 w_1)}{\partial x_1}
= \frac{\partial (w_1^{-1})}{\partial x_1} + w_1^{-1}\frac{\partial (x_1 w_1)}{\partial x_1}
=
$$
$$
=-w_1^{-1} \frac{\partial  w_1}{\partial x_1} + w_1^{-1} \left(1 +
x_1 \frac{\partial w_1}{\partial x_1} \right)=
w_1^{-1}\left( -\frac{\partial w_1}{\partial x_1} + 1 + x_1 \frac{\partial w_1}{\partial
x_1}\right) =
$$
$$
= w_1^{-1}\left(1 + (x_1 - 1) \frac{\partial w_1}{\partial
x_1}\right).
$$
Using the homomorphism $\varphi$ and taking into
account the following formulas
$$
(w_1^{-1})^{\varphi} = 1, ~~~\left(\frac{\partial w_1}{\partial
x_1}\right)^{\varphi} = 0,
$$
we obtain that
$$
\left(\frac{\partial (x_1 w)}{\partial
x_1}\right)^{\varphi} = 1.
$$
Then
$$
\frac{\partial (x_1 w)}{\partial x_2} = \frac{\partial (w_1^{-1} x_1 w_1)}{\partial x_2}
= \frac{\partial (w_1^{-1})}{\partial x_2} + w_1^{-1}\frac{\partial (x_1 w_1)}{\partial x_2}
=
$$
$$
=-w_1^{-1} \frac{\partial  w_1}{\partial x_2} + w_1^{-1} x_1 \frac{\partial w_1}{\partial x_2}=
-w_1^{-1} (1 - x_1) \frac{\partial w_1}{\partial
x_1}.
$$
Since  $\displaystyle \left(\frac{\partial w_1}{\partial
x_2}\right)^{\varphi} = 0,$ then $\displaystyle \left(\frac{\partial (x_1 w)}{\partial
x_2}\right)^{\varphi} = 0$.  Similarly,
$$
\frac{\partial (x_2 w)}{\partial x_1} = -w_1^{-1} (1 - x_2) \frac{\partial w_1}{\partial
x_1},~~~
\frac{\partial (x_2 w)}{\partial x_2} = w_1^{-1} \left( 1 + (x_2 - 1)
\frac{\partial w_1}{\partial
x_2}\right).
$$
Using the  homomorphism $\varphi$ we see that
$$
\widehat{\rho}_G (w(\varepsilon_{12}, \varepsilon_{21})) = 1.
$$
The proof is completed.

Now Lemma 4 and Lemma 6 imply Theorem 2. It then
follows that the representations $\rho :
\operatorname{IA}(F_n) \longrightarrow {\rm GL}_n(R)$
and $\widehat{\rho}_B : C_n \longrightarrow {\rm
GL}_n(\mathbb{Z}[t^{\pm 1}])$ are not faithful for all
$n\geq 2$.  Note that the question about
faithfulness of the Gassner representation of $P_n$,
$\geq 4$, is still open.

It is interesting to note that $Cb_2$ is generated by
elements $\varepsilon_{12}$ and $\varepsilon_{21}$ and
corresponding  matrices
$\widehat{\rho}_G(\varepsilon_{12}^{-1})$ and
$\widehat{\rho}_G(\varepsilon_{21}^{-1})$
are
$$
\widehat{\rho}_G(\varepsilon_{12}^{-1}) =
\left(
\begin{array}{cc}
\displaystyle \frac{\partial (x_1 \varepsilon_{12}^{-1})}{\partial x_1} &
\displaystyle \frac{\partial (x_1 \varepsilon_{12}^{-1})}{\partial x_2}  \\
& \\
\displaystyle \frac{\partial (x_2 \varepsilon_{12}^{-1})}{\partial x_1} &
\displaystyle \frac{\partial (x_2 \varepsilon_{12}^{-1})}{\partial x_2}   \\
\end{array}
\right)^{\varphi} =
\left(
\begin{array}{cc}
t_2 & 1 - t_1  \\
0 &  1    \\
\end{array}
\right),
$$
$$
\widehat{\rho}_G(\varepsilon_{21}^{-1}) =
\left(
\begin{array}{cc}
\displaystyle \frac{\partial (x_1 \varepsilon_{21}^{-1})}{\partial x_1} &
\displaystyle \frac{\partial (x_1 \varepsilon_{21}^{-1})}{\partial x_2}  \\
& \\
\displaystyle \frac{\partial (x_2 \varepsilon_{21}^{-1})}{\partial x_1} &
\displaystyle \frac{\partial (x_2 \varepsilon_{21}^{-1})}{\partial x_2}   \\
\end{array}
\right)^{\varphi}=
\left(
\begin{array}{cc}
1 &  0    \\
1 - t_2 & t_1  \\
\end{array}
\right).
$$
As we saw in the proof of the last Lemma these
matrices does not generate a free group.  But the
matrices
$$
\left(
\begin{array}{cc}
1 &  \lambda    \\
0 & 1  \\
\end{array}
\right),~~~
\left(
\begin{array}{cc}
1 &  0    \\
\mu & 1  \\
\end{array}
\right),~~~
\lambda, \mu \in \mathbb{C},~~~|\lambda | \geq 3,  |\mu| \geq 3,
$$
generate a group which is isomorphic to the free group
 $F_2$.



 \end{document}